\documentclass[psamsfont]{amsart}
\usepackage[utf8]{inputenc}
\usepackage{graphicx}
\usepackage[dvips]{epsfig}
\usepackage{pinlabel}
\usepackage{amsmath}
\usepackage{amsfonts}
\usepackage{latexsym}
\usepackage{amssymb}
\usepackage[usenames,dvipsnames]{color}
\usepackage{amsthm}
\usepackage[all]{xypic}
\usepackage{enumitem}
\usepackage[breaklinks=true]{hyperref}
\usepackage{tikz-cd}
\usepackage{xcolor}

\input xy

\xyoption{all}

\author{Emmy Murphy}
\address{Northwestern University, United States.}
\email{e\_murphy@math.northwestern.edu}

\theoremstyle{plain}
\newtheorem{thm}{Theorem}[section]
\newtheorem{lemma}[thm]{Lemma}
\newtheorem{prop}[thm]{Proposition}

\theoremstyle{definition}

\newtheorem{definition}[thm]{Definition}
\newtheorem{remark}[thm]{Remark}

\theoremstyle{remark}

\numberwithin{equation}{section}

\newcommand{\R}{\mathbb{R}}

\newcommand{\C}{\mathbb{C}}

\newcommand{\PP}{\mathbb{P}}

\newcommand{\J}{\mathcal{J}}
\newcommand{\FF}{\mathcal{F}}
\newcommand{\p}{\varphi}
\newcommand{\e}{\varepsilon}
\newcommand{\dd}{\partial}

\newcommand{\es}{\varnothing}
\newcommand{\op}{\operatorname}
\newcommand{\sse}{\subseteq}

\newcommand{\x}{\times}
\newcommand{\sm}{\setminus}
\newcommand{\pt}{\text{point}}

\newcommand{\wt}{\widetilde}

\newcommand{\ol}{\overline}

\newcommand{\std}{{\operatorname{std}}}

\newcommand{\Sh}{\operatorname{Sh}}

\newcommand{\Mor}{{\operatorname{Mor}}}

\newcommand{\Mod}{{\mathcal{M}od}}
\def\Op{{\mathcal O}{\it p}\,}

\begin{document}

\title{Arboreal singularities and loose Legendrians I}

\begin{abstract}
Arboreal singularities are an important class of Lagrangian singularities. They are conical, meaning that they can be understood by studying their links, which are singular Legendrian spaces in $S^{2n-1}_\std$. Loose Legendrians are a class of Legendrian spaces which satisfy an $h$--principle, meaning that their geometric classification is in bijective correspondence with their topological types. For the particular case of the linear arboreal singularities, we show that constructable sheaves suffice to detect whether any closed set of an arboreal link is loose.
\end{abstract}

\maketitle

\section{Introduction}\label{sec:intro}

Arboreal singularities, defined by Nadler \cite{N1}, are a class of Lagrangian singularities $L_T \sse (\R^{2n}, \omega_\std)$, corresponding to any rooted tree $T$ with no more than $n+1$ vertices. We review their definition in Section \ref{sec:defs}. $L_T$ is a conical singularity, meaning that if we define $\Lambda_T = L_T \cap S^{2n-1}$, then $\Lambda_T$ is a Legendrian complex in the contact sphere $(S^{2n-1}, \xi_\std)$, and $L_T$ is the cone of $\Lambda_T$, along the radial direction in $\R^{2n}_\std$.

Loose Legendrians, first defined in \cite{loose}, are notable because they satisfy an $h$--principle. This means that two loose Legendrians are isotopic among Legendrians whenever they are smoothly isotopic in a manner preserving the natural framings. Said differently, they contain certain local models which can be used to ``untangle'' any interesting geometry, and thus only their topology remains. They only exist in dimension $2n-1 \geq 5$: we will always assume this dimensional restriction throughout the paper. We review the definition and properties of loose Legendrians in Section \ref{sec:defs}.

Given contemporary tools, it is fairly easy to detect when a given Legendrian is non-loose: as soon as any Fukaya categorical invariant is non-vanishing -- such as the category of constructable sheaves or the Legendrian contact homology -- the Legendrian in question cannot be loose. Conversely, detecting when a Legendrian is loose seems to be an extremely difficult problem: all known conditions implying looseness are essentially reformulations of the definition. Throughout this paper we will work with constructable sheaf theory, as defined in \cite{STZ} building on work \cite{KS}. For any Legendrian $\Lambda \sse \R^{2n-1}_\std$, We denote by $\Sh(\Lambda)$ the derived category of constructable sheaves on $\R^n$ which are compactly supported and have singular support on $\Lambda$ (over any coefficient ring). Here $\R^n$ is identified as the front projection of $\R^{2n-1}_\std$.

Links of arboreal singularities $\Lambda_T$ are never loose: $\Sh(\Lambda_T)$ is equivalent to the category of derived modules of $T$ (thought of as a quiver), as shown by Nadler \cite{N1}. In particular there are many non-constant sheaves and it follows that $\Lambda_T$ is not loose. 

For a smooth Legendrian manifold, asking questions about the geometry of any non-trivial subset is not an interesting question. The reason is that any nontrivial closed subset of a manifold can be isotoped into a neighborhood of a space with positive codimension, and therefore the $h$-principle for subcritical isotropic embeddings gives a complete classification in correspondence with the smooth topology. However, since an arboreal link $\Lambda_T$ is a Legendrian \emph{complex} it has many interesting subsets which have nontrivial homology in the top dimension. This paper concerns the class of Legendrians obtained by taking any closed set of $\Lambda_{A_{n+1}}$, where $T = A_{n+1}$ is the linear tree.

\begin{thm}\label{thm:intro}
Let $\Lambda \sse \Lambda_{A_{n+1}}$ be any closed set. Then $\Lambda$ is loose if and only if $\Sh(\Lambda) \cong 0$.
\end{thm}

Here, a Legendrian complex is loose if every top dimensional cell is loose in the complement of all other cells. This is the natural definition of looseness for Legendrian complexes, as it ensures that the $h$-principle classification results apply. More generally, we will say that a given cell in a Legendrian complex is loose if that cell has a loose chart which is disjoint from all other cells.

The theorem is phrased in terms of closed sets in $\Lambda_{A_{n+1}}$, but similar to the case of smooth Legendrian manifolds the topology of these sets are mostly irrelevant. The only data they carry in terms of contact geometry is which cells of $\Lambda_{A_{n+1}}$ intersect $\Lambda$ in a proper set. Thus we will need an effective way to label these cells. Let $Q = A_{n+2}$ be the linear tree, thought of an appending a new zero object $0 \in Q$ to the tree $A_{n+1}$. We think of $Q$ as a quiver, meaning the category with $n+2$ elements and $\Mor(x, y)$ consisting of a single element if $x \leq y$ and $\Mor(x,y) = \es$ otherwise.

We claim then that there is a natural correspondence between the top-dimensional cells of $\Lambda_{A_{n+1}}$ and non-identity elements of $\Mor(Q)$, this is proven in Lemma \ref{lem:combin}. Thus, if $W \sse \Mor(Q)$, we can define a Legendrian $\Lambda_{Q[W^{-1}]}$ by deleting an open ball from any top-dimensional cell of $\Lambda_{A_{n+1}}$ which corresponds to (non-identity) elements of $W$. The following proposition follows quickly from the $h$-principle for subcritical isotropics.

\begin{prop}\label{prop:punct}
Let $\Lambda \sse \Lambda_{A_{n+1}}$ be a closed set. We define $W \sse \Mor(Q)$ as follows: for each $f \in \Mor(Q)$, $f \in W$ if and only if the top-dimensional cell of $\Lambda_{A_{n+1}}$ corresponding to $f$ intersects $\Lambda$ in a proper subset. Then $\Sh(\Lambda) = \Sh(\Lambda_{Q[W^{-1}]})$. Any given cell of $\Lambda$ is loose (rel $\Lambda$, see Definition \ref{def:loose rel}) if and only if the corresponding cell of $\Lambda_{Q[W^{-1}]}$ is loose (rel $\Lambda_{Q[W^{-1}]}$).
\end{prop}

The proof of Theorem \ref{thm:intro} then follows from the following three results, which are mostly independent. First we will prove a result generalizing the result from \cite{N1}. Denote by $\Mod(Q[W^{-1}])$ the derived category of modules $\rho: Q \to \op{Ch}_*$ sending $0 \in Q$ to $0 \in \op{Ch}_*$ and sending all morphisms in $W$ to quasi-isomorphisms. (i.e. $\Mod(Q[W^{-1}])$ is the category of derived modules of the localized category $Q[W^{-1}]$, preserving the initial object $0$.)

\begin{prop}\label{prop:shmod}
The category $\Sh(\Lambda_{Q[W^{-1}]})$ is equivalent to $\Mod(Q[W^{-1}])$.
\end{prop}

The next two results will both concern the notion of 2-out-of-6 closure, which will be an important concept in the paper. Given any subset of morphisms $W \sse Q$, we say that $W$ satisfies the \emph{2-out-of-6 property} if it contains all identities, and whenever we have a composition $a \overset{f}{\to} b \overset{g}{\to} c \overset{h}{\to} d$ so that $gf \in W$ and $hg \in W$, then $f$, $g$, $h$, and $hgf$ are all in $W$. We note that the 2-out-of-6 property implies the weaker 2-out-of-3 property: given a composition $a \overset{f}{\to} b \overset{g}{\to} c$, then whenever any two of the morphisms $f$, $g$, and $gf$ are in $W$, the remaining one is as well (we see this by inserting an identity). In particular if $W$ satisfies the 2-out-of-6 property then $W$ is closed under composition.

For any subset of morphisms $W \sse \Mor(Q)$, we denote by $\ol W \sse \Mor(Q)$ the \emph{2-out-of-6 closure} of $W$, i.e. the smallest subset of morphisms containing $W$ which satisfies the 2-out-of-6 property. The second main result of the paper is purely algebraic.

\begin{prop}\label{prop:loc quiv}
Let $Q = A_{n+2}$ be the linear quiver with initial object $0$ and let $W \sse \Mor(Q)$. For any $f \in \Mor(Q)$, $f \in \ol W$ if and only if for every module $\rho: Q[W^{-1}] \to \op{Ch}_*$, $\rho(f)$ is a quasi-isomorphism.
\end{prop}

The proof of Proposition \ref{prop:loc quiv} follows by constructing an explicit model of $Q[W^{-1}]$ (and Yoneda's Lemma). Finally, the remaining ingredient is to relate the above results to loose Legendrians.

\begin{prop}\label{prop:loose}
Let $D \sse \Lambda_{Q[W^{-1}]}$ be a top-dimensional cell, and let $f_D \in \Mor(Q)$ be the corresponding morphism. Then if $f_D \in \ol W$, it follows that $D$ is loose.
\end{prop}

The converse of the proposition is also true, as follows immediately from \cite{GKS} and \cite{loose}. Together these four propositions prove Theorem \ref{thm:intro}. In fact they prove a stronger result, by which we can work with each cell individually.

\begin{thm}\label{thm:main}
Let $\Lambda \sse \Lambda_{A_{n+1}}$ be any closed set, and let $D \sse \Lambda_{A_{n+1}}$ be any top-dimensional cell. If the inclusion functor $\Sh(\Lambda \sm D) \to \Sh(\Lambda)$ is an equivalence, it follows that $D$ is loose rel $\Lambda$.
\end{thm}

\subsection*{Structure of the paper}

In Section \ref{sec:defs} we review the necessary background for the paper, particularly the definitions and basic properties of arboreal singularities, loose Legendrians, and constructible sheaves. The following three sections prove the four main propositions above. 

In Section \ref{sec:prune} we define $\Lambda_{Q[W^{-1}]}$, and prove its basic properties described in Propositions \ref{prop:punct} and \ref{prop:shmod}. Section \ref{sec:algebra} contains the proof of Proposition \ref{prop:loc quiv}. Finally, Section \ref{sec:loose} contains the proof of Proposition \ref{prop:loose}.

\subsection*{Acknowledgments}

The author is grateful to the American Institute of Mathematics, for hosting a workshop on arboreal singularities in March 2018, and to D.~ Alvarez-Gavela, Y.~ Eliashberg, D.~ Nadler, and L.~ Starkston for stimulating discussions.

\section{Background}\label{sec:defs}

Throughout the paper we will always work with the contact manifold $\R^{2n-1}_\std$, whose contact structure is defined by the kernel of the $1$-form $dz - \sum_{i=1}^{n-1}y_idx_i$. This paper is concerned with Legendrian spaces in $\R^{2n-1}_\std$ which are more general than smooth manifolds. While they have a natural cellular structure, the geometry of the codimension $\geq 1$ portion is not interesting, and so we will use a naive definition which highlights the top-dimensional cells.

\begin{definition}\label{def:Legcomp}
A \emph{Legendrian complex} in a contact manifold $(Y, \xi)$ is a subset $\Lambda \sse Y$ which can be written as $\Lambda = \bigcup_i \Lambda_i$, where $\{\Lambda_i\}$ is a collection of smoothly embedded, connected Legendrian submanifolds with boundary and corners $\Lambda_i \sse (Y, \xi)$. $\{\Lambda_i\}$ are required to be mutually disjoint on their interiors. The set $\bigcup_i (\dd\Lambda_i)$ is called the \emph{singular set} of $\Lambda$.
\end{definition}

There are three basic ingredients necessary for the background of the paper: arboreal singularities, constructable sheaves, and loose Legendrians. For the latter two topics our treatment here is essentially standard (taken from \cite{STZ} and \cite{loose} respectively). Our treatment of arboreal singularities is somewhat novel, in that we define $\Lambda_T$ via its generic front projection, instead of the standard method of using conormals to hyperplanes in $\R^n$ \cite{N1}. This makes them appear perhaps less natural, but the advantage is that their contact geometry is more explicit.

\subsection{Arboreal singularities}\label{ssec:arbor}

Arboreal singularities are a class of Lagrangian singularities which have recently gained interest as a important class to understand, particularly within the context of skeleta of Weinstein manifolds. See \cite{N1, N2, Sta} for some important applications. We give a definition here which serves as a model for the links of these singularities, which are Legendrian complexes inside $\R^{2n-1}_\std = \dd B^{2n}_\std \sm \{\pt\}$.

For a fixed $n$, let $\Delta \sse \R^{n-1}$ be the standard embedding of the $(n-1)$-dimensional simplex, so that all $n$ vertices of $\Delta$ are equidistant from each other, and each vertex is distance $1$ away from the origin. The \emph{Venn diagram} is a configuration of $n$ round copies of $B^{n-1} \sse \R^{n-1}$, whose centers are the vertices of $\Delta$, all with equal radius $1+\e$ for a small $\e > 0$. Thus a given $k$--dimensional face of $\Delta$ contains a vertex $v$ in its closure if and only if the centroid of that face is contained in the ball corresponding to $v$. We define $r_v: \R^{n-1} \to [0, \infty)$ as the radial distance away from the point $v \in \R^{n-1}$. We also choose a bump function $\chi: [0, \infty) \to [0,1]$ which is non-increasing everywhere, equal to $r \mapsto (1+\e - r)^2$ for $r \in [1, 1+\e]$, equal to $0$ for $r \in [1+\e, \infty)$, and constant near $r=0$.  

Let $\pi: \R^{2n-1}_\std \to \R^n = \{(x_1, \ldots, x_{n-1}, z)\}$ be the front projection and let $p: \R^n \to \R^{n-1}$ be the projection onto the $x_i$ coordinates. Let $T$ be a rooted tree with $n+1$ vertices. We choose a bijection between the vertices of $\Delta$ with the non-root vertices of $T$. The \emph{arboreal link} corresponding to $T$ is a Legendrian $\Lambda_T \sse \R^{2n-1}_\std$, which is homeomorphic to a union of $S^{n-1}$ and $n$ copies of $D^{n-1}$, indexed by vertices $v \in T$. We define $\Lambda_T$ by defining its front projection.

\begin{figure}[h!]
\centering
\includegraphics[scale=0.6]{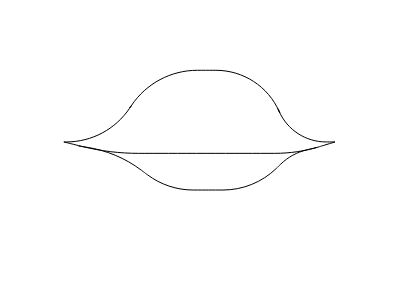}
\caption{The front of the standard Legendrian unknot $\pi(S^{n-1})$.}
\label{fig:unknot}
\end{figure}

The root $v_0$ of $T$ corresponds to $S^{n-1} \sse \Lambda_T$, where $\pi(S^{n-1})$ is the standard ``flying saucer'' front for the Legendrian unknot, see Figure \ref{fig:unknot}. We choose $\pi(S^{n-1})$ so that the lower branch of $\pi(S^{n-1})$ coincides with a large ball in the plane $\{z=0\}$, enough so that $p\pi(S^{n-1})$ contains the entire Venn diagram. We also choose the upper branch to have large $z$-value, in particular over the Venn diagram the upper branch should satisfy $z > n$.

For all other $v \in T$ which are not the root, the disk $D^{n-1}_v \sse \Lambda_T$ is defined so that $p \circ \pi: D_v^{n-1} \to \R^{n-1}$ takes (the interior of) $D_v^{n-1}$ diffeomorphically onto the ball in the Venn diagram centered at $v$. For all $v$, the interior of $\pi(D^{n-1}_v)$ will be contained in the open bounded region of $\R^n \sm \pi(S^{n-1})$, and the boundary of $\pi(D^{n-1}_v)$ only intersects the lower branch of $\pi(S^{n-1})$.

\begin{figure}[h!]
\centering
\includegraphics[scale=0.6]{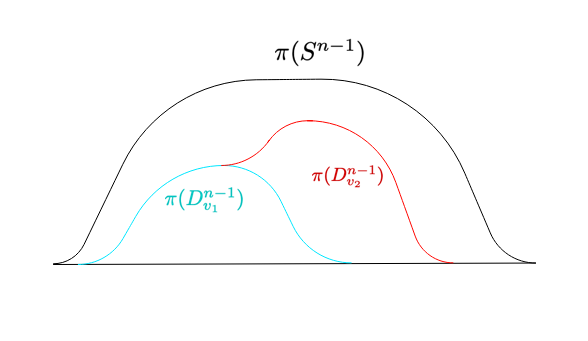}
\caption{The front of the standard $\Lambda_{A_3}$, the linear tree with three vertices.}
\label{fig:stack}
\end{figure}

\begin{figure}[h!]
\centering
\includegraphics[scale=0.6]{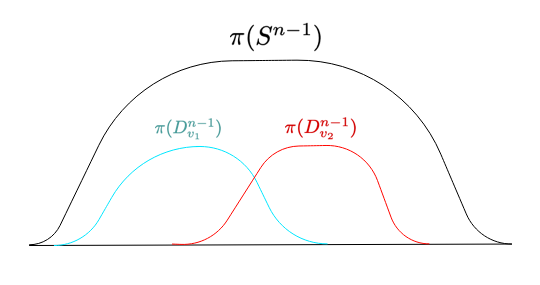}
\caption{The front of the standard $\Lambda_T$ where $T$ is the unique non-linear rooted tree with three vertices (i.e. the linear tree except with the root being the center vertex.)}
\label{fig:cross}
\end{figure}

We complete the definition inductively with respect to the partial ordering given by $T$. We define $\pi(D^{n-1}_v) = \{z = \sum_{w \leq v} \chi(r_w), r_v \in [0,1+\e]\}$. Informally, each $\pi(D^{n-1}_v)$ is a ``dome'' which is stacked on top of all previous domes sitting below it in $T$. Given two vertices $v_1$ and $v_2$ which are incomparable in $T$, $\pi(D^{n-1}_{v_1})$ intersects $\pi(D^{n-1}_{v_2})$ in a $(n-2)$ disk, but they are transverse on their interiors and therefore the interiors of $D^{n-1}_{v_1}$ and $D^{n-1}_{v_2}$ are disjoint in $\R^{2n-1}_\std$. See Figures \ref{fig:stack} and \ref{fig:cross}. It would be instructive to the reader to make sure they can clearly picture all $\pi(\Lambda_T) \sse \R^3$, corresponding to the $4$ distinct rooted trees with $n+1 = 4$ vertices.

\begin{remark}\label{rmk:technical}
We make a number of technical remarks that may be useful to experts in the theory, but do not have any bearing on the main result.

It is perhaps not immediately obvious that this definition is equivalent to Nadler's original definition from \cite{N1}: we claim that, if $L_T \sse B^{2n}_\std$ is Nadler's arboreal singularity corresponding to the rooted tree $T$, then $L_T \cap \dd B^{2n}_\std = \Lambda_T$, under the identification $\R^{2n-1}_\std \cong S^{2n-1}_\std \sm \{\pt\}$. We will not prove this fact here, but we give an intuitive sketch. The boundary of the standard $\R^n \sse \R^{2n}$ is isotopic to the standard Legendrian unknot, and the coordinate hyperplanes in $\R^n$ intersect the sphere at infinity along the Venn diagram in $S^{n-1}$. Similar to the original definition, the front projection is by definition a conormal construction, and so it suffices to work entirely in $\R^n$, but this correspondence involves a stereographic projection. To make the details precise, one would have to keep track on signs well (to see that ``positive conormal'' in the original definition corresponds to the positive conormal implicit in the front projection), and pay detailed attention to the singularities along the regions where $\Lambda_T$ is glued from disks.


Since their original appearance, arboreal singularities have been generalized to signed versions, see \cite{Sta}. Our definition here corresponds to arboreal singularities which are purely positive. One can give Legendrian definitions similar to the above section for arboreal links with arbitrary signs on the edge set, though we will not do this here.

$\chi$ has a discontinuous second derivative at the point $1+\e$, thus the Legendrian curve with front $\{z = \chi(x)\}$ is continuous, but not $C^1$ smooth at $x = 1+\e$. Therefore, as an instance of our definition above, if $D^2_v$ is defined by the front $\{z = \chi(r_v) + \chi(r_w)\}$, the closed disk $D^2_v$ is not smooth: it has a corner at its boundary point $r_v = r_w = 1+\e$. More generally, $D^{n-1}_v$ with have order $k$ corners whenever there are $k$ non-root vertices $w$ satisfying $w \leq v$. Regardless of $n$ and $k$, the front projection $\pi(D^{n-1}_v)$ is a disk which is $C^1$ smooth but not $C^2$ everywhere on its boundary, and $p: \pi(D^{n-1}_v) \to \R^{n-1}$ is a $C^1$ diffeomorphism onto a closed round ball. Both $\pi$ and $p \circ \pi$ are $C^\infty$ diffeomorphisms on the interior of $D^{n-1}_v$.


As a Legendrian complex, $\Lambda_T$ has many singularities itself. The singularities of $\Lambda_T$ correspond to the intersection points of the spheres in the Venn diagram. In particular $\Lambda_T$ has $2n$ isolated singularities: the isolated intersections in the Venn diagram are $(n-1)$--fold, that is to say they involve every sphere except one, and every subset of $n-1$ spheres in the Venn diagram intersects in two points. If we take the Lagrangian projection $\R^{2n-1}_\std \to \R^{2n-2}_\std$, each isolated singularity of $\Lambda_T$ gives a isolated Lagrangian singularity: we claim that this is the arboreal singularity corresponding to the subtree $T \sm \{v\}$ where $v$ is the vertex in $T$ corresponding to the unique disk $D^{n-1}_v$ whose boundary does not intersect the singularity.\footnote{Here ``subtree'' is in the sense of full subcategories of the quiver, rather than subgraphs of a directed graph: in a rooted tree, any collection of vertices containing the root gives rise to a subtree, which is itself a rooted tree.} In particular, even within the self-contained definitions here we could ask whether the Legendrian link of the Lagrangian projection of a singularity of $\Lambda_T$ is a lower dimensional arboreal link. We will also neglect to prove this.

If $T$ is a rooted tree with $k < n+1$ vertices, we can still define a Legendrian $\Lambda_T \sse \R^{2n-1}_\std$, simply by attaching $k$ copies of $D^{n-1}$ to $S^{n-1}$ along a $k$--component subdiagram of the Venn diagram. As above, we note without proof that $\Lambda_T$ is the link of the arboreal singularity $L_T \x \R^{n-k} \sse \R^{2k}_\std \x \R^{2n-2k}_\std$.
\end{remark}

\subsection{Constructable sheaves in contact geometry}\label{ssec:sheaves}

In this section we review the basic material concerning constructable sheaves in contact geometry. Our account here principally follows \cite{STZ}.

We fix a ring $R$. Throughout this paper a \emph{sheaf} $\FF$ on a manifold $M$ will always refer to a chain complex of sheaves of $R$-modules. Such a sheaf is called \emph{constructable} if it is locally constant with respect to some stratification $\mathcal{S}$ of $M$. For such a sheaf $\FF$, we will define the singular support of $\FF$ in terms of stratified Morse theory. If $x \in M$ is contained in a neighborhood $U$ and $f: U \to \R$ is a smooth function, we define the \emph{Morse group} of $(x, f)$ to be the cone of the restriction map $\FF(U \cap f^{-1}(-\infty, \delta)) \to \FF(U \cap f^{-1}(-\infty, -\delta))$, where $\delta > 0$ is a small positive number. Implicit in this definition is the fact that, for sufficiently small $U$ and sufficiently small $\delta$ (allowed to depend on $U$), this chain complex does not depend on $U$ or $\delta$.

If $\FF$ is locally constant with respect to the stratification $\mathcal{S}$ and $f: U \to \R$ is a function which is stratified Morse with respect to $\mathcal{S}$, then whenever the Morse group of $(x, f)$ is not acyclic we say that the covector $df_x \in T^*M$ is \emph{characteristic}. Thus a covector $p \in T^*M_x$ is characteristic if there exists a function $f$ with $df_x = p$ which is stratified Morse with respect to $\mathcal{S}$ and whose Morse group is cohomologically non-trivial. This definition does depend on the choice of $\mathcal{S}$, but the \emph{singular support} --- defined to be the point-set closure of all characteristic covectors --- only depends on $\FF$. We denote this set by $SS(\FF) \sse T^*M$. Then $SS(\FF)$ is Lagrangian everywhere it is smooth \cite{KS}, furthermore it is conical with respect to fiberwise radial dilation, i.e. the standard Liouville vector field. Thus if $S^*M$ is the sphere bundle of $T^*M$ (at infinity or with a chosen metric), $SS(\FF) \sse S^*M$ is a Legendrian space. If $\pi: S^*M \to M$ is the front projection and $\mathcal{S}$ is any stratification of $M$ for which $\FF$ is locally constant, then $\pi(SS(\FF))$ is contained in the codimension $\geq 1$ strata of $\mathcal{S}$.

\begin{definition}\label{def:sheaves}
Let $\Lambda \sse S^*M$ be a Legendrian complex, then we define the category $\Sh(\Lambda)$ as follows. Objects $\FF \in \Sh(\Lambda)$ are chain complexes of sheaves of $R$-modules, which are constructable and satisfy $SS(\FF) \sse \Lambda$. We also require that $\FF$ is cohomologically bounded at each stalk (i.e. perfect), and that it has compact support in the case where $M$ is non-compact.

Morphisms in $\Sh(\FF)$ consist of derived morphisms, that is usual morphisms of sheaves, where additionally all quasi-isomorphisms are localized.

In the particular case where $\Lambda \sse \R^{2n-1}_\std$ (which is the only case we consider), we use the canonical embedding $\R^{2n-1}_\std \cong S^*_-\R^n \sse S^*\R^n$ to define $\Sh(\Lambda)$. Here $S^*_-\R^n$ consists of those ($\R_+$ projectivized) covectors which evaluate negatively on the vector field $\dd_z$ on $\R^n$. We note that the front projection $\pi: \R^{2n-1}_\std \to \R^n$ is equal to the base projection $\pi: S^*_-\R^n \to \R^n$ in this correspondence.
\end{definition}

An important fact, proved in \cite{GKS}, is that $\Sh(\Lambda)$ is a Legendrian invariant: it only depends on $\Lambda$ up to contact isotopy. 

\begin{thm}[\cite{GKS}]\label{thm:GKS}
Any contact isotopy $\p_t: S^*M \to S^*M$ induces an equivalence of categories $\Sh(\Lambda) \to \Sh(\p_1(\Lambda))$
\end{thm}

This theorem is not directly relevant to the results in this paper, but it is significant in order to quickly prove the converses to the main theorems here: the converse of Theorem \ref{thm:main} and the ``only if'' portion of Theorem \ref{thm:intro} both follow immediately from Theorems \ref{thm:GKS} and \ref{thm:close loose}.

\begin{remark}
For smooth Legendrians $\Lambda$, any family of Legendrian embeddings is induced by an ambient contact isotopy. For Legendrian complexes this is false. This will not be very relevant for us, since we are always working with Legendrians which are subsets of a fixed Legendrian complex $\Lambda_{A_{n+1}}$. In more generality, asking whether the singular sets of Legendrian complexes are contact isotopic essentially reduces to the question of whether they have contactomorphic neighborhoods: since a neighborhood of the singular set is contained in a small neighborhood of a subcritical isotropic space, $h$-principles can promote local contactomorphisms to global contact isotopies of the singular sets. Extending these isotopies to the interior of the top-dimensional cells then proceeds as in the smooth case.
\end{remark}

If $\mathcal{S}$ is a stratification of $\R^n$ which refines $\pi(\Lambda)$, then any sheaf in $\Sh(\Lambda)$ is locally constant with respect to $\mathcal{S}$ (up to quasi-isomorphism). Typically it is easy to arrange that $\mathcal{S}$ consists of finitely many cells, and that each cell is contractible. Let $\mathcal{Q}_{\mathcal{S}}$ be the finite thin category defined by the combinatorics of $\mathcal{S}$: the objects of $\mathcal{Q}_{\mathcal{S}}$ are the cells of $\mathcal{S}$, and each $\Mor(C_1, C_2)$ contains either one or zero elements, according to whether $C_1$ is in the closure of $C_2$. Let $\Mod(\mathcal{Q}_{\mathcal{S}})$ denote the category of derived perfect $R$-modules of $\mathcal{Q}_{\mathcal{S}}$. Since any sheaf in $\Sh(\Lambda)$ is constant on each cell of $\mathcal{S}$, we have a fully faithful embedding $\Sh(\Lambda) \to \Mod(\mathcal{Q}_{\mathcal{S}})$.

We will refrain from stating any general result here (and we discard the notation from the previous paragraph), but the intuitive principle is useful to keep in mind: \emph{$\Sh(\Lambda)$ is equivalent to a full subcategory of $\Mod(\mathcal{Q})$, for some finite thin category $\mathcal{Q}$.} For Legendrian knots, this is done in detail in Sections 3.3, 3.4, and 3.5 of \cite{STZ}. In more generality figuring out which $\mathcal{Q}$ and which full subcategory of $\Mod(\mathcal{Q})$ correctly calculates $\Sh(\Lambda)$ involves analyzing the singularities of $\pi(\Lambda)$. We do this for our relevant cases in Section \ref{sec:prune}.

\subsection{Loose Legendrians}\label{ssec:loose}

Among Legendrian submanifolds $\Lambda \sse (Y, \xi)$ of dimension $\dim Y = 2n-1 \geq 5$ there exists the class of loose Legendrians. There are two important points to make about loose Legendrians:
\begin{itemize}
\item[--]They are defined in terms of containing a model. That is, there is a universal Legendrian $\Lambda_\ell \sse \R^{2n-1}_\std$, equal to the standard plane $\{z=y=0\}$ outside of a compact set, so that $\Lambda$ is loose if and only if there exists an open set $U \sse Y$, so that the pair $(U, U \cap \Lambda)$ is contactomorphic to $(\R^{2n-1}_\std, \Lambda_\ell)$. Thus, while looseness is a \emph{global} property (Legendrians have no local invariants), it is semi-local in the sense that it can be \emph{certified} by exhibiting it is loose on a single open set. 
\item[--]Loose Legendrians are classified up to Legendrian isotopy by data which is purely diffeo-topological. They are typically thought of as being geometrically trivial: Fukaya categories and constructable sheaf categories are all trivial for loose Legendrians. They are useful particularly for constructions: if you want to show that a loose Legendrian $\Lambda$ satisfies Property $\mathfrak{X}$ it suffices to build any Legendrian $\wt\Lambda$ which is loose and satisfies Property $\mathfrak{X}$. As long as the mild topological constraints are satisfied then $\Lambda$ will be isotopic to $\wt\Lambda$, and therefore $\Lambda$ also satisfies Property $\mathfrak{X}$ as long as the property is invariant up to contact isotopy.
\end{itemize}

In order to explain these informal descriptions we will need the notion of a \emph{formal Legendrian isotopy}. Recall that when $(Y, \xi)$ is any contact manifold, the vector bundle $\xi$ (forgetting its embedding $\xi \sse TY$) is equipped with a linear symplectic form, which is well-defined up to conformal scaling, when $\xi = \ker \alpha$ this symplectic form is $d\alpha|_{\xi}$. Since $\alpha|_\Lambda = 0$ implies $d\alpha|_\Lambda = 0$, it follows that any Legendrian, simply by virture of being tangent to $\xi$ everywhere, must in fact be a Lagrangian subspace $T\Lambda_x \sse \xi_x$ for all $x \in \Lambda$.

\begin{definition}
A \emph{formal Legendrian embedding} is a pair $(f, F_s)$, where $f: \Lambda\to Y$ is a smooth embedding, and $F_s: T\Lambda \to TY$ is a homotopy of bundle monomorphisms covering $f$ for all $s \in [0,1]$. $F_s$ is required to connect $F_0 = df$ and $F_1$ -- a map whose image $F_1(T\Lambda)$ lies inside $\xi$ as a Lagrangian subspace.

A Legendrian embedding is precisely a formal Legendrian embedding which is constant in $s$. In particular a \emph{formal Legendrian isotopy} between Legendrian embeddings $f_0, f_1: \Lambda \to (Y, \xi)$ is a path in the space of formal Legendrian embeddings connecting $f_0$ to $f_1$.
\end{definition} 

We remark that the forgetful map from formal Legendrian embeddings to smooth embeddings is a Serre fibration. Furthermore the homotopy fiber can be identified with something explicit, such as $\op{Map}(\Omega\Lambda, O_{2n+1}/U_n)$ for the stably parallelizable case (the typical case is the gauge group of the Lagrangian Grassmannian of the symplectic bundle $T^*\Lambda \otimes \C \to \Lambda$).

As alluded to above, a Legendrian $\Lambda \sse (Y, \xi)$ is called \emph{loose} if there is an open $U \sse Y$ so that the pair $(U, U \cap \Lambda)$ is contactomorphic to $(\R^{2n-1}_\std, \Lambda_\ell)$, where the Legendrian $\Lambda_\ell \sse \R^{2n-1}_\std$ is a standard model called a \emph{loose chart}. The specific geometry of the model $\Lambda_\ell$ will be relevant to us soon, we define it below at Proposition \ref{prop:zzloose}. One basic property of $\Lambda_\ell$ is that $\dim(\Lambda_\ell) \geq 2$, therefore \emph{by definition a loose Legendrian must sit inside a contact manifold $(Y, \xi)$ of dimension $2n-1 \geq 5$.} The following theorem is the main result from \cite{loose}.

\begin{thm}\label{thm:basic loose}
Let $f_0, f_1: \Lambda \to (Y, \xi)$ be two Legendrian embeddings of a connected smooth manifold $\Lambda$, and assume they are formally Legendrian isotopic. If they are both loose, then they are Legendrian isotopic.
\end{thm}

Through the rest of the section we explain how to generalize Theorem \ref{thm:basic loose} to Legendrians which are not connected, smooth manifolds. 

\begin{definition}\label{def:loose rel}
Let $\Lambda$ be a smooth connected Legendrian and let $A \sse Y$ is some closed set, possibly intersecting $\Lambda$. We say that $\Lambda$ is \emph{loose rel $A$} if there is a loose chart $U \sse Y$ for $\Lambda$ so that $U \cap A \sse U \cap \Lambda$.

If $\Lambda = \bigcup_i \Lambda_i$ is a Legendrian complex, we say that $\Lambda$ is \emph{loose} if each $\Lambda_i$ is loose rel $\Lambda$.
\end{definition}

For a simple example, if $\Lambda_i$ are disjointly embedded smooth submanifolds, then $\Lambda$ is just a Legendrian link with components $\{\Lambda_i\}$. Then to say that the link $\Lambda$ is loose means that each component is loose, with a loose chart disjoint from the other components. We note that it is fairly easy to construct a non-loose link, whose components are individually loose.

We say that two Legendrian complexes $\Lambda, \wt \Lambda \sse (Y, \xi)$ are \emph{formally Legendrian isotopic} if, firstly, there is an ambient contact isotopy $\p_t: Y \to Y$ sending an open neighborhood of the singular set of $\Lambda$ to a neighborhood of the singular set of $\wt \Lambda$, and secondly that there is a formal isotopy supported compactly on the interior of all $\p_1(\Lambda_i)$ which take $\p_1(\Lambda_i)$ to $\wt\Lambda_i$. 

\begin{thm}\label{thm:complex loose}
Let $\Lambda, \wt\Lambda \sse (Y, \xi)$ be two Legendrian complexes which are formally Legendrian isotopic. If they are both loose, then they are ambient contact isotopic.
\end{thm}

The proof of Theorem \ref{thm:complex loose} just follows by applying the theorem to each $\Lambda_i$ individually. We state a more technical version of the loose Legendrian classification to make this clear.

\begin{thm}[\cite{loose}]\label{thm:close loose}
Let $f_0, f_1: \Lambda \to (Y, \xi)$ be two Legendrian embeddings of a connected manifold $\Lambda$, which are equal on a closed set $A \sse Y$. Suppose further that $f_0 = f_1$ on an open set $U \sse Y$, $U \cap A = \es$, and that $(U, U \cap f_0(\Lambda)) = (U, U \cap f_1(\Lambda))$ is a loose chart. Suppose that there is a formal Legendrian isotopy $(g_t, g_{s, t})$ between $f_0$ and $f_1$, which is supported on $Y \sm (A \cup U)$. We assume that every connected component of $\Lambda \sm A$ intersects $U$.

Then there is a Legendrian isotopy $f_t: \Lambda \to Y$ connecting $f_0$ to $f_1$, with the following properties. $f_t$ is supported on $Y \sm A$, and outside of $U$ $f_t$ is $C^0$ close to $g_t$.
\end{thm}

From this it is clear that we can work cell by cell to prove Theorem \ref{thm:complex loose}. Though we did not assume that the loose charts of $\Lambda_i$ were necessarily equal for $\Lambda$ and $\wt \Lambda$, this can be arranged since all Darboux balls with smooth boundary are isotopic, via an isotopy supported away from any Legendrian disjoint from it.

One important aspect of loose Legendrian complexes is that the loose charts \emph{themselves} do not have to be disjoint. If $\Lambda = \Lambda_1 \cup \Lambda_2$ is a Legendrian complex and there are loose charts $U_1, U_2 \sse Y$ with $U_1 \cap \Lambda_2 = \es = U_2 \cap \Lambda_1$, and $(U_1, U_1 \cap \Lambda_1) \cong (\R^{2n-1}_\std, \Lambda_\ell) \cong (U_2, U_2 \cap \Lambda_2)$, then this implies that $\Lambda$ is loose in the strongest sense. Even if it is not assumed that $U_1 \cap U_2 = \es$, we can guarantee that there are other loose charts $U_i' \sse Y$ for $\Lambda_i$ satisfying the same properties as $U_i$ above (for $i=1, 2$), with the additional property that $U_1' \cap U_2' = \es$.

Here is the proof: let $\wt U_2 \sse Y$ be any set which is disjoint from $\Lambda_1$, $U_1$, and $U_2$, and which intersects $\Lambda_2$ in the set $\{z=y_i=0\} \sse \R^{2n-1}_\std$, i.e. the standard zero section in $\J^1\R^{n-1}$. Any small neighborhood of a point is contactomorphic to such a set so this is immediate. Let $\wt\Lambda_2 \sse Y$ be the Legendrian which agrees with $\Lambda_2$ outside of $\wt U_2$, and inside $\wt U_2$ we replace the zero section $\{z=y_i=0\}$ with the compactly supported loose Legendrian (formally representing the zero section). Notice that $\Lambda_2$ and $\wt\Lambda_2$ are both loose Legendrians and they are formally isotopic, therefore they are isotopic.

So there is a contact flow $\p_t: Y \to Y$ with $\p_1(\Lambda_2) = \wt\Lambda_2$. In fact since $U_2$ is a loose chart for both $\Lambda_2$ and $\wt\Lambda_2$, we can assume that $\p_t$ is supported outside of the set $V = U_2 \cup \Op(\gamma) \cup \wt U_2$, where here $\Op(\gamma)$ is the neighborhood of an arc in the interior of $\Lambda_2$, chosen so that $V$ is connected. This follows since the formal isotopy between $\Lambda_2$ and $\wt \Lambda_2$ happens in $\wt U_2$. Notice also that $\p_t$ is supported outside of a neighborhood of $\Lambda_1$. Then let $U_2' = \p_1^{-1}(\wt U_2)$, and let $U_1' = \p_1^{-1}(U_1)$. 

\begin{prop}\label{prop:freebloose}
Let $\Lambda = \bigcup_i\Lambda_i$ be a Legendrian complex, and let $\Lambda_k \sse \Lambda$ be a cell with free boundary: i.e. there is an open set $U$ around a point $x \in \dd \Lambda_k$ so that $U \cap \Lambda_i = \es$ for all $i \neq k$. Then $\Lambda_k$ is loose rel $\Lambda$.
\end{prop}

\begin{proof}
Let $U$ be a small Darboux coordinate ball intersecting $\Lambda_k$ is the zero section and disjoint from all other $\Lambda_i$, and let $\wt \Lambda_k$ be the Legendrian obtained from $\Lambda_k$ by replacing the portion of $\Lambda_k$ inside $U$ with a Legendrian which is loose and formally isotopic to the zero section. Certainly $\wt \Lambda$ is loose rel $A$, thus it suffices to show that $\Lambda_k$ is Legendrian isotopic (rel $A$) to $\wt \Lambda_k$.

$\Lambda_k$ and $\wt \Lambda_k$ are formally Legendrian isotopy, this implies in particular that they are formally regular homotopic: i.e. there is a family of maps $f_t: \Lambda_k$ interpolating between the inclusion, and the embedding of $\wt \Lambda$, and furthermore these maps are covered by injective bundle maps $F_t: T\Lambda_k \to TY$ with Lagrangian image. Thus, the $h$-principle for Legendrian immersions \cite{hbook} implies that there is a family of Legendrian immersions $g_t: \Lambda_k \to Y$, fixed on $A$, interpolating between the inclusion and the embedding of $\wt \Lambda_k$.

If $g_t$ is generically perturbed then all double points of the maps $g_t: \Lambda_k \to Y$ occur at isolated times $t \in [0,1]$, at isolated points in $\Lambda_k$. Let $x_1, \ldots x_m \in \Lambda_k$ be a list of these points (ignoring the times the double points occur). There is a smooth isotopy through inclusions $\p_t: \Lambda_k \to \Lambda_k$ which is fixed on $A \cap \Lambda_k$, and so that $x_1, \ldots x_m \notin \p_1(\Lambda_k)$. Here we use the fact that $\dd \Lambda_k \sm A \neq \es$.

Any smooth isotopy of inclusions $\p_t: \Lambda_k \to \Lambda_k$ is realized by a ambient contact isotopy, acting on any contact manifold where $\Lambda_k$ lies as a Legendrian. Thus we can find a contact isotopy $\p_t: Y \to Y$, so that $\p_1(\Lambda) \sse \Lambda$ contains none of the points $x_1, \ldots, x_m$. Define a Legendrian isotopy as follows: first let the contact isotopy $\p_t$ act from $t \in [0,1]$, then we concatenate with the family $g_t\circ \p_1: \Lambda_k \to Y$ which consist of Legendrian embeddings, then finally concatenate with the reverse flow of $\p_t$.
\end{proof}

Notice that this proposition shows that any Legendrian with free boundary satisfies an $h$-principle, \emph{even rel boundary}. For instance, if we take the unit disk $D^{n-1}_\std \sse \R^{n-1} = \{z=p=0\} \sse \R^{2n-1}_\std$, Proposition \ref{prop:freebloose} implies that there is a loose chart $U \sse \R^{2n-1}_\std$for $D^{n-1}_\std$ with $U \cap \dd D^{n-1}_\std = \es$. Then the $h$-principle for loose Legendrians, applied rel $\dd D^{n-1}_\std$, classifies all Legendrians equal to $D^{n-1}_\std$ near a neighborhood of the boundary. This cannot be done for example for $D^{n-1}_\std \sse \J^1D^{n-1}$, because this Legendrian is not loose.

Finally, we will need a concrete definition of a loose chart, in order to understand the relationship to arboreal singularities in later sections. There are many possible definitions, we will use one that is useful for our purpose.

\begin{figure}[h!]
\centering
\includegraphics[scale=0.6]{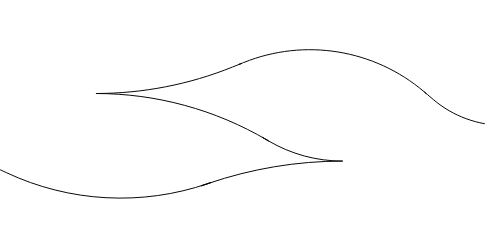}
\caption{The Legendrian zig-zag $\wt \Lambda_0 \sse \R^3_\std$.}
\label{fig:zz}
\end{figure}

\begin{prop}[\cite{loose, affine}]\label{prop:zzloose}
Let $\wt \Lambda_0\sse \R^3_\std$ be the standard Legendrian zig-zag as given in Figure \ref{fig:zz}. Consider the Legendrian $\Lambda_\ell = \wt \Lambda_0 \x D^{n-2} \sse \R^3_\std \x T^*D^{n-2}$. Then $(U, \Lambda_\ell)$ is loose.
\end{prop}

\begin{remark}
Though Proposition \ref{prop:zzloose} is not the original definition of a loose chart, it can be taken as the definition: the pair $(U, \Lambda_\ell)$ contains a loose chart (according to the original definition of loose chart), and also any loose chart contains a contactomorphic copy of $(U, \Lambda_\ell)$. Since looseness of any Legendrian is defined by containing a contactomorphic copy of a loose chart, it follows that a Legendrian is loose if and only if it contains a contactomorphic embedding of the pair $(U, \Lambda_\ell)$. In this sense it is reasonable to \emph{define} a loose chart to be the pair $(U, \Lambda_\ell)$. This is the perspective we will take throughout the paper.

However, while Proposition \ref{prop:zzloose} will be important for us as our one touchstone of a specific Legendrian which is loose, it will not be important to us that a loose chart contains the model $(U, \Lambda_\ell)$. Ultimately it is less important what the model for a loose chart \emph{is}, compared to what it \emph{does} (i.e. Theorem \ref{thm:close loose}). 
\end{remark}

\section{Pruning}\label{sec:prune}

We now focus on arboreal links $\Lambda_T$ in the particular case where $T = A_{n+1}$ is linear. First we analyze the combinatorial structure of $\Lambda_{A_{n+1}}$. Throughout we work with open cells $C = \op{Int}(C) \sse \Lambda_{A_{n+1}}$.

As defined in Section \ref{ssec:arbor}, the front of the arboreal link $\Lambda_{A_{n+1}}$ is defined as a union of disks $D^{n-1}_v = \{z = \sum_{w \leq v}\chi(r_w)\}$. In particular, whenever $v_1 > v_2$ in the tree $T$, we see that in the front projection the interior of $D^{n-1}_{v_1}$ is disjoint from $D^{n-1}_{v_2}$, since the $z$ coordinate is strictly larger everywhere. Since $T = A_{n+1}$ is linear, this implies that all disks $D^{n-1}_v$ have disjoint interiors.

\begin{figure}[h!]
\centering
\includegraphics[scale=0.6]{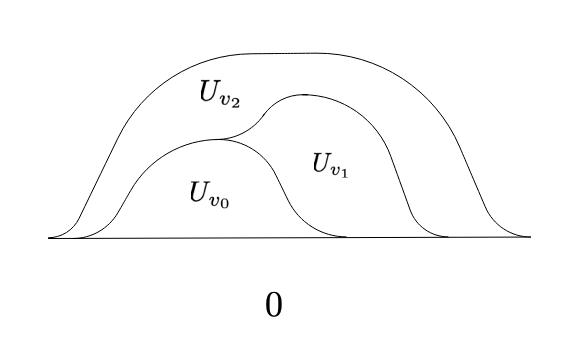}
\caption{The open sets $U_j \sse \R^n$. This example corresponds to $A_3 = \{v_0 \to v_1 \to v_2\}$; recall that $v_0$, the root of $A_3$ is distinct from $0$, the initial object of $Q$.}
\label{fig:regions}
\end{figure}

Thus we see that $\R^n \sm \pi(\Lambda_{A_{n+1}})$ has exactly $n+1$ bounded components: the initial flying saucer $\pi(S^{n-1}) \sse \pi(\Lambda_{A_{n+1}})$ has one bounded component in its complement, and each disk $D^{n-1}_v$ divides a single component in two. More concretely, since each $\pi(D^{n-1}_v)$ is attached on \emph{top} of $\pi(S^{n-1})$ and the previous $\pi(D^{n-1}_w)$, there is a unique component $U \sse \R^n \sm \pi(\Lambda_{A_{n+1}})$ lying below $\pi(D^{n-1}_v)$ and containing it in its closure. We denote this component by $U_{v-1}$, where $v-1$ is the vertex of $A_{n+1}$ preceding $v$. In particular, the root $v_0 \in A_{n+1}$ is associated to the bounded component $U_{v_0}$ lying below every $\pi(D^{n-1}_w)$. If $v \in A_{n+1}$ is the maximal element, we define $U_v$ to be the component of $\R^n \sm \pi(\Lambda_{A_{n+1}})$ lying above each $\pi(D^{n-1}_w)$, i.e. the unique component whose closure intersects the upper hemisphere of $\pi(S^{n-1})$. See Figure \ref{fig:regions}.

The unbounded component of $\R^n \sm \pi(\Lambda_{A_{n+1}})$ will also be important for our purposes. Thus we define $Q$ to be the quiver obtained by appending an initial object $0$ to $A_{n+1}$. Thus $Q = A_{n+2}$, though we avoid this notation to avoid confusion with the Legendrian in one larger dimension. We define $U_0 \sse \R^n \sm \pi(\Lambda_{A_{n+1}})$ to be the unbounded component, and thus the components of $\R^n \sm \pi(\Lambda_{A_{n+1}})$ are in correspondence with the vertices of $Q$. We claim that the top dimensional cells of $\Lambda_{A_{n+1}}$ are naturally in bijective correspondence with all morphisms of $Q$. More generally:

\begin{lemma}\label{lem:combin}
Let $2 \leq k \leq n+1$, and let $\{v_1, \ldots, v_k\}$ be any collection of vertices in $Q$. Then $\ol U_{v_1} \cap \ldots \cap \ol U_{v_k}$ is the image of the closure of a single cell of $\Lambda_{A_{n+1}}$, whose codimension in $\Lambda_{A_{n+1}}$ is $k-1$. (Here $\ol U_v$ denotes the point-set closure of $U_v$.) Every cell of $\Lambda_{A_{n+1}}$ arises uniquely in this way, and thus the $m$-cells of $\Lambda_{A_{n+1}}$ are in correspondence with $(n-m+1)$-element subsets of vertices of $Q$.
\end{lemma}

\begin{proof}
If $C \sse \Lambda_{A_{n+1}}$ is any $m$-cell, we can define $V_C = \{v \in Q; \pi(C) \sse \ol U_v \}$. Then it is clear that $\ol \pi(C) \sse \bigcap_{v \in V_C} \ol U_v$. Thus to prove the entire lemma it suffices to show the reverse inclusion.

We have either $C \sse \op{Int}(D^{n-1}_w)$ for some $w$, or else $C \sse S^{n-1} \sse \Lambda_{A_{n+1}}$. In the former case $\pi(D^{n-1}_w) \sse \ol U_{w-1}$, and $w-1 \in V_C$ is the minimal element in $V_C$, by the definition of $U_{w-1}$. In the latter case $\pi(S^{n-1}) \sse \ol  U_0$, and again $0 \in V_C$ is the minimal element of $V_C$. Since $\ol U_{w-1} \cap \ol U_v \sse \pi(D^{n-1}_w)$ for any $v > w$, it follows that $\bigcap_{v \in V_C} \ol U_v \sse \pi(D^{n-1}_w)$. Again addressing the latter case we have that $\ol U_0 \cap \pi(\Lambda_{A_{n+1}}) = \pi(S^{n-1})$, and so $\bigcap_{v \in V_C} \ol U_v \sse \pi(S^{n-1})$.

Let $D = D^{n-1}_w$ or $D = S^{n-1}$, according to the cases above so that $C \sse D$. Let $x \in D \sse \Lambda_{A_{n+1}}$ be any point not in $\ol C$. Then there is a $u \in A_{n+1}$ (not the root) so that $\dd D^{n-1}_u \cap D$ separates $C$ from $x$. The two components of $\pi(D \sm \dd D^{n-1}_u)$ are $\ol U_{u-1} \cap D$ and $D \sm \ol U_{u-1}$. Thus $\pi(C) \sse \ol U_{u-1}$ if and only if $\pi(x) \notin \ol U_{u-1}$, and so it follows that $\pi(x) \notin \bigcap_{v \in V_C} \ol U_v$.
\end{proof}

The next lemma shows that the above identification is natural with respect to the isomorphism $\Sh(\Lambda_{A_{n+1}}) \cong \Mod(Q)$ from \cite{N1}. Recall that $\Mod(Q)$ is defined to be the derived category of modules sending $0 \in Q$ to $0 \in \op{Ch}_*$.

\begin{lemma}\label{lem:sheaves}
There is an equivalence of categories $\mathcal{N}: \Sh(\Lambda_{A_{n+1}}) \to \Mod(Q)$ satisfying the following property. Let $V \sse Q$ be any full subcategory and let $U_V = \op{Int}(\bigcup_{v \in V} \ol U_v)$. Then for any two $\FF_1, \FF_2 \in  \Sh(\Lambda_{A_{n+1}})$, $\FF_1|_{U_V} \cong \FF_2|_{U_V}$ if and only if $\mathcal{N}(\FF_1)|_V \cong \mathcal{N}(\FF_2)|_V$.
\end{lemma}

\begin{proof}
We start by defining the functor $\mathcal{N}$, so let $\FF \in \Sh(\Lambda)$. For each object $v \in Q$, we define $\mathcal{N}(\FF)(v)$ to be the chain complex $\FF(U_v)$. For each morphism $\p: v_1 \to v_2$, choose a point $x_\p$ in the interior of the cell $\ol U_{v_1} \cap \ol U_{v_2}$. Let $U_\p$ be a small neighborhood of the poiont $x_\p$, and consider the restriction maps $\FF(U_\p) \to \FF(U_\p \cap U_{v_1})$ and $\FF(U_\p) \to \FF(U_\p \cap U_{v_2})$. The former map is necessarily a quasi-isomorphism: its cone is the Morse group of the pair $(x_\p, f)$ for some function $f: U_\p \to \R$ which is positive on $U_\p \cap U_{v_2}$ and negative on $U_\p \cap U_{v_1}$. Since $U_{v_1}$ lies below $U_{v_2}$ this means that $df_{x_\p}$ evaluates positively on $\dd_z$, and thus the point $(x_\p, \PP(df_{x_\p})) \in S^*\R^n$ is not contained in $\R^{2n-1}_\std = S^*_-\R^n$. In particular this point is not contained in $\Lambda_{A_{n+1}}$, so it is not in the singular support of $\FF$, so $\FF(U_\p) \to \FF(U_\p \cap U_{v_1})$ is a quasi-isomorphism.

Furthermore since $\pi(\Lambda_{A_{n+1}})$ is disjoint from $U_{v_1}$ and $U_{v_2}$, the restriction maps $\FF(U_{v_1}) \to \FF(U_\p \cap U_{v_1})$ and $\FF(U_{v_2}) \to \FF(U_\p \cap U_{v_2})$ are quasi-isomorphisms. Thus we have a map (in the derived category) 
$$\FF(U_{v_1}) \to \FF(U_\p \cap U_{v_1}) \to \FF(U_\p) \to \FF(U_\p \cap U_{v_2}) \to \FF(U_{v_2}),$$
where the second and fourth maps are inverted quasi-isomorphisms. This defines $\mathcal{N}(\FF)(\p) \in \Mor(\Mod(Q))$. To see that $\mathcal{N}(\FF)$ respects composition we choose morphisms $\p: v_1 \to v_2$ and $\psi: v_2 \to v_3$ in $Q$, and choose a point $x$ in the interior of the cell $\ol U_{v_1} \cap \ol U_{v_2} \cap \ol U_{v_3}$. A neighborhood $U_x$ of $x$ in $\R^n$ is contactomorphic to the standard trivalent Legendrian front $\pi(\Lambda_3) \sse \R^2$, extended trivially by $\R^{n-2}$. See Figure \ref{fig:vertex}.

\begin{figure}[h!]
\centering
\includegraphics[scale=0.45]{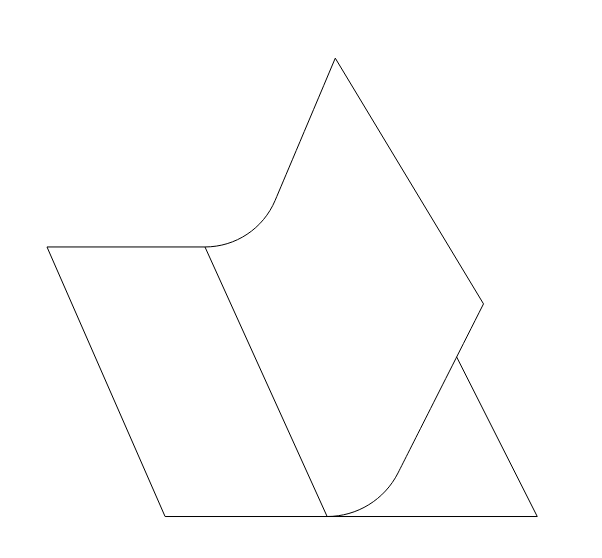}
\caption{A trivalent vertex, trivially extended.} 
\label{fig:vertex}
\end{figure}

The analysis of $\FF$ near the point $x_{\p, \psi}$ is similar to the analysis of the cusp singularity done in \cite[Section 3.3]{STZ}. The combinatorics of the restriction maps is expressed in the commutative diagram:
$$\begin{tikzcd}
 & U_{v_3} & U_\psi \arrow{l} \arrow{d} \\
 &  & U_{v_2} \\
U_{\psi \circ \p}  \arrow{dr}  \arrow{uur} & U_x \arrow{uu}  \arrow{uur} \arrow{ur} \arrow{l} \arrow{r} \arrow{d} & U_\p \arrow{u} \arrow{dl} \\
 & U_{v_1} & 
\end{tikzcd}$$
The three downward arrows $\FF(U_\psi) \to \FF(U_{v_2})$, $\FF(U_\p) \to \FF(U_{v_1})$ and $\FF(U_{\psi \circ \p}) \to \FF(U_{v_1})$ must be quasi-isomorphisms, this is the same as the argument above. Since $\mathcal{N}(\FF)(\psi \circ \p)$ is the morphism $\FF(U_{v_1}) \to \FF(U_{v_3})$ obtained by following the left hand side of the diagram, and $\mathcal{N}(\FF)(\psi) \circ \mathcal{N}(\FF)(\p)$ is obtained by following the right hand side, we see that this maps are equal and so $\mathcal{N}(\FF)$ is a module.

$\mathcal{N}$ is a functor in the obvious way: a map between sheaves $\FF_1 \to \FF_2$ consists in particular of maps of chain complexes $\FF_1(U_v) \to \FF_2(U_v)$ for all $v$, since these maps commute with all restrictions this defines a map between modules of $Q$. $\mathcal{N}$ is a faithful functor since every stalk of $\FF$ is naturally quasi-isomorphic to $\FF(U_v)$ for some $v$ (this holds for any Legendrian in $S^*_-\R^n$).

Let $M \in \Mod(Q)$. We show that there is a sheaf $\FF_M \in \Sh(\Lambda_{A_{n+1}})$ with $\mathcal{N}(\FF_M) = M$. Besides showing that $\mathcal{N}$ is essentially surjective, the construction will be natural and thus show that $\mathcal{N}$ is a full functor. It suffices to construct $\FF_M$ in the case where $M = \mathcal{Y}^v = R^{\Mor(v, \cdot)}$ is the Yoneda module of $v$, a chosen vertex.\footnote{Here we are using the fact that $\Sh(\Lambda)$ consists of perfect sheaves, and $\Mod(Q)$ consists of perfect modules. In fact in this case the situation is so explicit that this is not necessary: arbitrary $\mathcal{C}$--valued sheaves with singular support in $\Lambda_{A_{n+1}}$ are equivalent to $\op{Func(Q, \mathcal{C})}$. But since perfect derived modules are the most relevant structure in symplectic geometry we take advantage of the structure, since it involves less analysis of the singularities of $\Lambda_{A_{n+1}}$.}

Let $U_{\geq v} = \bigcup_{w \geq v} \ol U_w$, and let $\wt U_{\geq v} \sse U_{\geq v}$ denote the complement of the ``bottom boundary'': a point $x \in \dd U_{\geq v}$ is also in $\wt U_{\geq v}$ unless there is a function $f$ defined near $x$ which is positive on $\op{Int}(U_{\geq v})$, negative outside $U_{\geq v}$, and satisfying $\frac{\dd f}{\dd z} > 0$. If $R$ is the constant sheaf on $\wt U_{\geq v}$ and $\FF = \iota_!(R)$, where $\iota: \wt U_{\geq v} \to \R^n$ is the inclusion map, then the singular support of $\FF_M$ is contained in $\Lambda_{A_{n+1}}$, and clearly $\mathcal{N}(\FF) = \mathcal{Y}^v$.

The statement that $\mathcal{N}$ is an equivalence when intertwined with restrictions $V \sse Q$ and $U_V \sse \R^n$ is obvious from the definition.


\end{proof}

Let $W \sse \Mor(Q)$ be an arbitrary set, which corresponds to some collection of top-dimensional cells as in Lemma \ref{lem:combin}. Recall that $\Lambda_{Q[W^{-1}]}$ is defined to be the Legendrian obtained from $\Lambda_{A_{n+1}}$ by deleting a small open ball from the interior of each cell corresponding to a morphism of $W$.

\begin{proof}[Proof of Proposition \ref{prop:shmod}:]
Since $\Lambda_{Q[W^{-1}]} \sse \Lambda_{A_{n+1}}$, we have the inclusion functor $\Sh(\Lambda_{Q[W^{-1}]}) \to \Sh(\Lambda_{A_{n+1}})$ which is fully faithful. For a morphism $f: v_1 \to v_2$ define the set $U_f = \op{Int}(\ol U_{v_1} \cup U_{v_2})$, so that pair $(U_f, \Lambda_{A_{n+1}})$ is diffeomorphic to $(\R^n, \R^{n-1} \x \{0\})$. If $f \in W$, then $U_f \sm \Lambda_{Q[W^{-1}]}$ is connected, and so any $\FF \in \Sh(\Lambda_{Q[W^{-1}]})$ must be locally constant on $U_f$. Thus $\mathcal{N}(\FF)(f)$ is a quasi-isomorphism, and so $\mathcal{N}: \Sh(\Lambda_{Q[W^{-1}]}) \to \Mod(Q)$ factors through $\Mod(Q[W^{-1}])$.

It remains to show that $\mathcal{N}: \Sh(\Lambda_{Q[W^{-1}]}) \to \Mod(Q[W^{-1}])$ is essentially surjective. If $\FF \in \Sh(\Lambda_{A_{n+1}}$ is any sheaf with $\mathcal{N}(\FF) \in \Mod(Q[W^{-1}])$, then for all $f \in W$ $f: v_1 \to v_2$, we have that $\mathcal{N}(\FF)|_{\{v_1, v_2\}}$ is quasi-isomorphic to a constant module, and therefore $\FF_{U_f}$ is quasi-isomorphic to a constant sheaf. Thus the singular support of $\FF$ is disjoint from $U_f$, and therefore $\FF \in \Sh(\Lambda_{Q[W^{-1}]})$ since $\Lambda_{Q[W^{-1}]} = \Lambda_{A_{n+1}}$ outside of such $U_f$.
\end{proof}

\begin{proof}[Proof of Proposition \ref{prop:punct}:]
Let $\Lambda \sse \Lambda_{A_{n+1}}$ be any closed set, and let $W \sse \Mor(Q)$ consist of those morphisms $f \in \Mor(Q)$ corresponding to those top-dimensional cells $D_f$ so that $D_f \cap (\Lambda_{A_{n+1}} \sm \Lambda) \neq \es$. 

First, assume that $\Lambda$ contains the entire singular set of $\Lambda_{A_{n+1}}$ it its interior, i.e. all cells of codimension at least $1$ are contained in $\op{Int}(\Lambda)$. Then the same proof as above shows that the fully faithful image of $\Sh(\Lambda)$ under $\mathcal{N}$ is exactly the subcategory $\Mod(Q[W^{-1}])$: the inclusion $\Lambda \sse \Lambda_{Q[W^{-1}]}$ shows that $\mathcal{N}(\Sh(\Lambda)) \sse \Mod(Q[W^{-1}])$ and the reverse inclusion follows since $\Lambda = \Lambda_{A_{n+1}}$ outside of $\bigcup_{f\in W}U_f$. The inclusion $\Lambda \sse \Lambda_{Q[W^{-1}]}$ shows that every cell which is loose rel $\Lambda_{Q[W^{-1}]}$ is also loose rel $\Lambda$. Since the interior of every top-dimensional cell is diffeomorphic to $D^{n-1}$, we can choose an ambient contact isotopy $\p_t$ so that $\p_t(\Lambda_{Q[W^{-1}]}) \sse \Lambda_{Q[W^{-1}]}$ and $\p_1(\Lambda_{Q[W^{-1}]}) \sse \Lambda$: $\p_t$ is the map which expands the punctures of $\Lambda_{Q[W^{-1}]}$ radially until the boundary of the punctures lie in a small neighborhood of the singular set. Then any cell which is loose rel $\Lambda$ is also loose rel $\p_1(\Lambda_{Q[W^{-1}]})$ and therefore also loose rel $\Lambda_{Q[W^{-1}]}$. This completes the proof under the stated assumption.

Finally, let $V \sse \Lambda_{A_{n+1}}$ be the closure of a small neighborhood of the singular set, and for arbitrary $\Lambda \sse \Lambda_{A_{n+1}}$ let $\Lambda' = \Lambda \cup V$. Then $V \sm \Lambda$ retracts, via contact isotopy fixed on $\Lambda$, to an arbitrarily small neighborhood of a subcritical isotropic complex. In particular any loose chart which is disjoint from $\Lambda$ can be made disjoint from $\Lambda'$, via contact isotopy fixed on $\Lambda$. Also, any sheaf whose singular support is contained in $\Lambda'$ in fact has singular support contained in $\Lambda$. Thus the assumption in the previous paragraph loses no generality.
\end{proof}

\section{Localization of quivers}\label{sec:algebra}

In this section we give a proof of Proposition \ref{prop:loc quiv}. The following is clearly a necessary condition.

\begin{lemma}\label{lem:loc quiv}
Let $W \sse \Mor(Q)$. Then for any $f \in W$, the image of $f$ under the functor $Q \to Q[W^{-1}]$ is an isomorphism only if $f \in \ol W$.
\end{lemma}

In fact the lemma immediately implies Proposition \ref{prop:loc quiv}: if $\rho(f)$ is a quasi-isomorphism for all $\rho: Q[W^{-1}] \to \op{Ch}_*$ the Yoneda lemma implies that $f \in Q[W^{-1}]$ is an isomorphism, by which we conclude $f \in \ol W$ from the lemma. The converse follows immediately from the definition of 2-out-of-6 closure: if in the composition $a \overset{f}{\to} b \overset{g}{\to} c \overset{h}{\to} d$ we know that $gf$ and $hg$ are isomorphisms, then $g^{-1} = f(gf)^{-1} = (hg)^{-1}h$ is a two-sided inverse for $g$, by which it immediately follows that $f$, $h$, and $hgf$ are isomorphisms as well. 

\begin{proof}[Proof of Lemma \ref{lem:loc quiv}:]
If $W$ admitted a calculus of fractions this would be a known result \cite[7.1.20]{KS}. Since it does not, in order to establish the result we will have to construct a model for $Q[W^{-1}]$. It is clear that $Q[W^{-1}] = Q[\ol W^{-1}]$, so without loss of generality we assume that $W = \ol W$.

This model is as follows: for any $a, b \in Q$, we define $\Mor(a, b)$ to consist of equivalence classes of diagrams $a \overset{f}{\to} m \overset{w}{\leftarrow} m_0 \overset{g}{\to} b$, where $w \in W$ and $f, g$ are arbitrary morphisms. We denote this morphism by the formal expression $gw^{-1}f$. The equivalence relation is defined by the two horizontal morphisms in the diagram being equivalent:

\begin{equation}\label{eq:equiv}
\begin{tikzcd}
a \arrow{r}{f} \arrow{d}{1_a} & m \arrow{d}{x} & m_0 \arrow{l}{w} \arrow{r}{g} \arrow{d}{x_0} & b \arrow{d}{1_b} \\
a \arrow{r}{f'} & m' & m_0' \arrow{l}{w'} \arrow{r}{g'} & b
\end{tikzcd}
\end{equation}

where $w, w' \in W$ and all other morphisms are arbitrary. That is, whenever $xw = w'x_0$, we have $(g'x_0)w^{-1}f = g(w')^{-1}(xf)$, and in particular $g1_m^{-1}f = g' 1_{m'}^{-1} f'$ whenever $g'f' = gf$.

Given morphisms $a \overset{f}{\to} m \overset{w}{\leftarrow} m_0 \overset{g}{\to} b$ and $b \overset{h}{\to} n \overset{u}{\leftarrow} n_0 \overset{k}{\to} c$ we need to define the composition $ku^{-1}h \circ gw^{-1}f \in \Mor(a, c)$. If $y_0 \in \Mor(m_0, n_0)$, then $hg = uy_0$ in $Q$ (since $\Mor(m_0, n)$ has at most one element), and therefore we define $ku^{-1}h \circ gw^{-1}f = (ky_0)w^{-1}f$. Similarly if $y \in \Mor(m, n)$ we define $ku^{-1}h \circ gw^{-1}f = ku^{-1}(yf)$. If both $y_0$ and $y$ exist, the definition is unambiguous because $uy_0 = yw \in \Mor(m_0, n)$. If neither exist, then since $Q$ is linear we must have both $\Mor(n_0, m_0)$ and $\Mor(n, m)$ are non-empty, and thus we have a diagram of the form $n_0 \overset{z_0}{\to} m_0 \overset{hg}{\to} n \overset{z}{\to} m$. Since $hgz_0 = u \in W$ and $zhg = w \in W$, and since $W$ satisfies the 2-out-of-6 property, it follows that $zhgz_0 \in W$ and we define $ku^{-1}h \circ gw^{-1}f = k(zhgz_0)^{-1}f$.

One easily sees that this composition is well defined on equivalence classes by doing a case analysis. 

This defines a category $Q[W^{-1}]$ with a natural functor $Q \to Q[W^{-1}]$ sending $f \in \Mor(a,b)$ to $f1_a^{-1}1_a = 1_b1_b^{-1}f$, and whenever $f \in W$ we have an inverse in $Q[W^{-1}]$ given by $1_b f^{-1} 1_a$. It remains to show that the only $f \in \Mor(Q)$ which are sent to isomorphisms are already in $W$ (this also establishes the universal property of the localization $Q[W^{-1}]$, though this is essentially obvious). 

First we make an auxiliary claim, that the morphism $a \overset{f}{\to} m \overset{w}{\leftarrow} m_0 \overset{g}{\to} a$ is equivalent to the identity only if $f, g \in W$. We show that this property is preserved under the equivalence given in Diagram \ref{eq:equiv} when $a=b$. First we suppose $z \in \Mor(m_0', m)$, then $m_0 \overset{x_0}{\to} m_0' \overset{z}{\to} m \overset{x}{\to} m'$ is a 2-out-of-6 diagram with $zx_0 = w \in W$ and $xz = w' \in W$, so $x, x_0 \in W$. Since $f' = xf$ and $g = g'x_0$, by the 2-out-of-3 property we see that if either $f, g \in W$ or $f', g' \in W$, then both hold. Instead, supposing $\Mor(m_0', m) = \es$ implies that $\Mor(m, m_0') \neq \es$ and so let $z' \in \Mor(m, m'_0)$. But then we have a diagram $a \overset{f}{\to} m \overset{z'}{\to} m_0' \overset{g'}{\to} a$, so by linearity of $Q$ we have $m = m_0' = a$ and $f = g' = 1_a$. Thus $g = w$ and $f' = w'$, so $f, g, f', g' \in W$.

Having established the claim, suppose that $f \in \Mor(a, b)$ and it has an inverse $ku^{-1}h$ in $\Mor(Q[W^{-1}])$, so $(fk)u^{-1}h = 1_b$ and $ku^{-1}(hf) = 1_a$. The claim implies that $fk, h, k, hf \in W$, and thus the 2-out-of-3 property ensures $f \in W$.
\end{proof}

\section{Criteria for looseness}\label{sec:loose}

\begin{figure}[h!]
\centering
\includegraphics[scale=0.6]{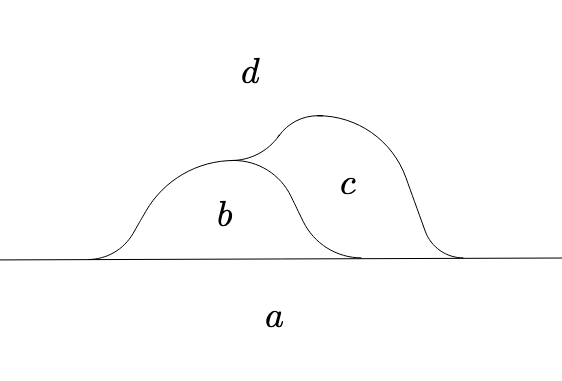}
\caption{The front $\pi(\Lambda_6)$.}
\label{fig:L6}
\end{figure}

In this section we present the proof of Proposition \ref{prop:loose}. First we prove a basic lemma which is a geometric model for the 2-out-of-6 property. We let $\Lambda_6 \sse \R^3_\std$ be the Legendrian $1$-complex given in Figure \ref{fig:L6}, this is essentially $\Lambda_{A_3}$ except in ``long knot'' format. Notice in particular that the category $\wt{\Sh}(\Lambda_6)$ consisting of constructable sheaves with singular support on $\Lambda_6$ and \emph{no conditions on the support} is equivalent to modules of the category $a \overset{f}{\to} b \overset{g}{\to} c \overset{h}{\to} d$. First, we show that any algebraic subcategory can be represented by a geometric slice. Recall that $\pi: \R^{2n-1}_\std \to \R^n$ is the front projection.

\begin{lemma}\label{lem:A4 slice}
Let $\mathcal{C} \sse Q$ be any full subcategory which is equivalent to the category $A_4 = a \overset{f}{\to} b \overset{g}{\to} c \overset{h}{\to} d$. Then there is a set $U \sse \R^n$, so that the pair $(\pi^{-1}(U), \pi^{-1}(U) \cap \Lambda_{A_{n+1}})$ is contactomorphic to $(\R^3_\std \x T^*D^{n-2}, \Lambda_6 \x D^{n-2})$, and so that the restriction of constructable sheaves from $\R^n$ to $U$ realizes the restriction of modules from $Q$ to $\mathcal{C}$.
\end{lemma}

\begin{proof}
This is essentially a corollary of Lemma \ref{lem:sheaves}, but to be explicit we will prove it more geometrically. Let $B^{n-1}_v \sse \R^{n-1}$ be the ball in the standard Venn diagram corresponding to $v$, where $v$ is any vertex in $A_{n+1}$ excepting the root. Thus for $v \in Q$, the ball $B^{n-1}_{v+1}$ is defined in the obvious way unless $v = 0$, or $v \in Q$ is the maximal element. In particular $B^{n-1}_{b+1}$ and $B^{n-1}_{c+1}$ definitely exist. 

We choose a compact arc $\gamma \sse \R^{n-1}$ as follows. Firstly we require $\gamma$ is completely disjoint from $B^{n-1}_{v+1}$, whenever $v$ is not equal to $a, b, c,$ or $d$. $\gamma$ will be completely contained in $B^{n-1}_{a+1} \cap B^{n-1}_{d+1}$ as long as both of these balls are defined, if not we require $\gamma$ to be contained in $B^{n-1}_{a+1}$ or $B^{n-1}_{d+1}$ if only one of these is defined (or else impose nothing if $a = 0$ and $d$ is the maximum). The first endpoint of $\gamma$ is required to be outside $B^{n-1}_{b+1}$ and $B^{n-1}_{c+1}$. As we follow $\gamma$, it is required to enter $B^{n-1}_{b+1}$, then enter $B^{n-1}_{c+1}$, then exit $B^{n-1}_{b+1}$, then exit $B^{n-1}_{c+1}$. See Figure \ref{fig:arc6}.

\begin{figure}[h!]
\centering
\includegraphics[scale=0.5]{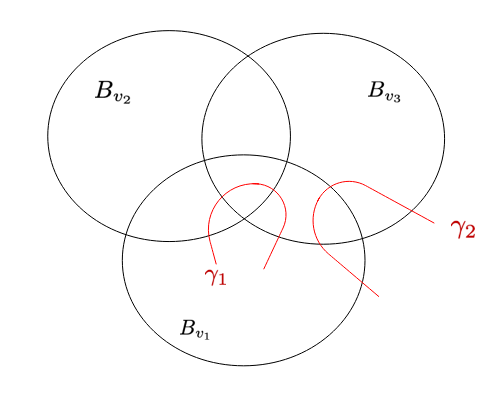}
\caption{Possible examples for the arc $\gamma$, pictured here in the case $n=3$. Here $\gamma_1$ is the example corresponding to the subcategory $v_0 \to v_1 \to v_2 \to v_3$, and $\gamma_2$ corresponds to the subcategory $0 \to v_0 \to v_2 \to v_3$. Notice that every such $\gamma$ is a open set of a small linking circle of a codimension $2$ singularity, as Lemma \ref{lem:sheaves} would predict.}
\label{fig:arc6}
\end{figure}

In $\R^n$, the $2$-plane $\gamma \x \R$ intersects $\pi(\Lambda_{A_{n+1}}$ transversely. Following from the definition of $\Lambda_{A_{n+1}}$, we see that the intersection is diffeomorphic as a pair to $\pi(\Lambda_6) \sse \R^2$, with some additional disjoint arcs lying completely above or below the diffeomorphic copy of $\pi(\Lambda_6)$. Namely there are additional intersections: a slice of the upper hemisphere of $\pi(S^{n-1})$, a slice of $\pi(D^{n-1}_{d+1})$ if $d$ is not the maximum, and a slice of the lower hemisphere of $\pi(S^{n-1})$, which is a trivial arc whenever $a \neq 0$ (if $a=0$ the lower hemisphere of $\pi(S^{n-1})$ forms a portion of our copy of $\Lambda_6$). Since all these additional arcs lie above or below the diagram, we can chose a compact piece of $2$-plane $P \sse \gamma \x \R$ so that the only intersection is the copy of $\pi(\Lambda_6)$

Since the intersection is transverse we have a tubular neighborhood $\wt P \sse \R^n$ of $P$ so that $\wt P \cap \pi(\Lambda_{A_{n+1}}) \cong \pi(\Lambda_6 \x B^{n-2})$. Then $U = \pi^{-1}(\wt P)$ is the desired neighborhood.
\end{proof}

Let $\Lambda_0 \sse \Lambda_6$ be the Legendrian $1$-complex obtained by puncturing the two edges corresponding to the morphisms $gf$ and $hg$.

\begin{lemma}\label{lem:loose slice}
$\Lambda_0 \x D^{n-2} \sse \R^3_\std \x T^*D^{n-2}$ is loose.
\end{lemma}

\begin{proof}
Throughout this proof all isotopies will be with compact support. If instead of $\Lambda_0$ we started with the standard zig-zag $\wt\Lambda_0 \sse \R^3$ (the front of a smooth Legendrian curve), this is Proposition \ref{prop:zzloose}. Let $\Lambda = \Lambda_0 \x D^{n-2}$ and $\wt \Lambda = \wt \Lambda_0 \x D^{n-2} \sse \R^3_\std \x T^*D^{n-2}$, which are equal outside of a neighborhood of the codimensional $1$ stratum of $\Lambda$. Let $\Lambda_\ell \sse \R^{2n-1}_\std$ be a (smooth) loose Legendrian which is formally isotopic to the standard zero section, and equal to it outside a compact set. Let $\wt\Lambda_1$ be the Legendrian built from $\wt \Lambda$ by implanting four small copies of $\Lambda_\ell$, near points corresponding to the top-dimensional cells of $\Lambda$. See Figure \ref{fig:supaloose}. Since both $\wt \Lambda$ and $\wt \Lambda_1$ are loose Legendrians, there is a Legendrian isotopy taking $\wt\Lambda$ to $\wt\Lambda_1$. Since the $h$-principle for loose Legendrians is for parametrized Legendrians, we can assume that the isotopy at time $1$ fixes $\wt\Lambda$ \emph{pointwise} everywhere it fixes it setwise, i.e. outside of the implanted copies of $\Lambda_\ell$.

\begin{figure}[h!]
\centering
\includegraphics[scale=0.22]{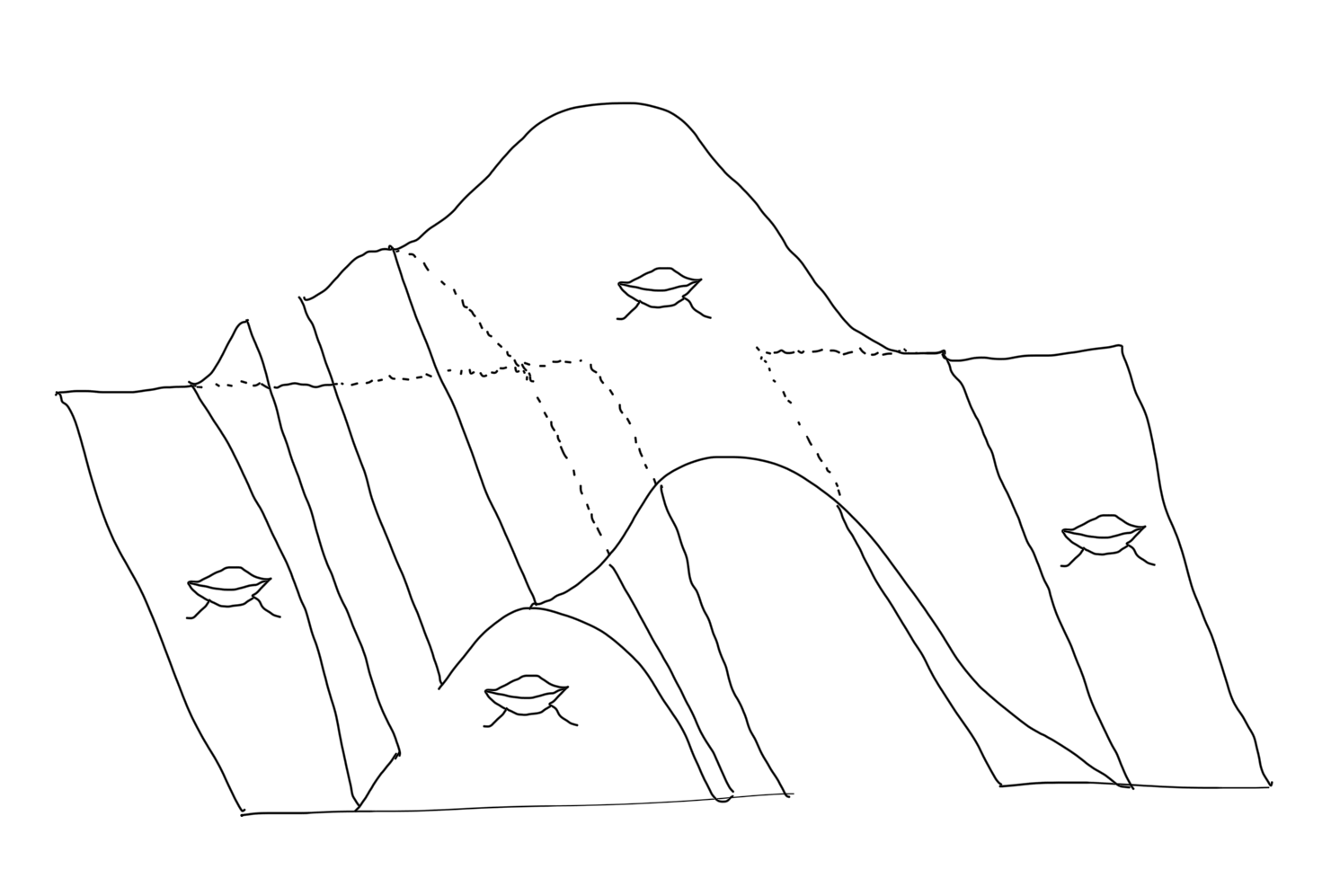}
\caption{The Legendrian $\p_1(\Lambda)$, which is explicitly loose. The isotopy $\p_t$ taking this Legendrian to $\Lambda$ is built using the fact that the large scale structure of $\Lambda_0$, is $C^0$ close to a zig-zag itself, as explained in the proof.}
\label{fig:supaloose}
\end{figure}

Let $\p_t: \R^3_\std \x T^*D^{n-2}$ be an extension of this isotopy to an ambient contact isotopy. $\Lambda = \wt \Lambda$ outside of an open set $U \sse \R^3_\std \x T^*D^{n-2}$, so that $\p_1$ is the identity on $U \cap \wt \Lambda$. By further ambient isotopy if needed, we can assume that every point of $\Lambda$ within the support of $\p_t$ is $\e$-close to a point of $\wt \Lambda$. Then, since $\p_1$ is continuous, we see that for every point $x \in \Lambda \cap U$, $\p_1(x)$ is arbitrarily close to $x$. In particular $\p_1(\Lambda) = \p_1(\wt \Lambda) = \wt\Lambda_1$ outside of a small open neighborhood of the codimension $1$ stratum. In particular the implanted copies of $\Lambda_\ell$ are embedded in $\p_1(\Lambda)$. Thus each of the four top-dimensional cells in $\p_1(\Lambda)$ has a loose chart, disjoint from the other cells. Since $\p_t$ is a contact isotopy these loose charts also exist in $\Lambda$.
\end{proof}

Finally we prove a lemma that allows us to set up an inductive proof.

\begin{lemma}\label{lem:loose induct}
Let $\Lambda$ be a Legendrian complex (in any contact manifold), and let $C(\Lambda)$ be the collection of top-dimensional cells of $\Lambda$. Let $C_0, C_1 \sse C(\Lambda)$ be two disjoint collections. Then, if every cell in $C_0$ is loose rel $\Lambda$, and every cell in $C_1$ is loose rel $\Lambda \sm C_0$, then every cell in $C_0 \cup C_1$ is loose rel $\Lambda$.
\end{lemma}

\begin{proof}
For each $D_i \in C_1$, let $U_i$ be a loose chart for $D_i$, which only intersects cells in $C_0$ (besides $D_i$ itself). Since $U_i$ is a ball with smooth boundary, there is no obstruction to finding a \emph{smooth} isotopy of the singular set of $\Lambda$, so that after the isotopy the singular set of $\Lambda$ is is disjoint from all $U_i$. The singular set is itself a subcritical isotropic complex, and the space of formal isotropic embeddings is a Serre fibration over the space of smooth embeddings, thus the $h$-principle for subcritical isotropics \cite{hbook} implies that there is a isotopy through isotropics disjoining the singular set of $\Lambda$ from $U_i$. We can then use the isotopy extension theorem to find a contact isotopy $Y \to Y$ which is fixed on $\Lambda \sm C_0$, so that after the isotopy each $U_i$ is disjoint from the singular set of $\Lambda$. Thus we can assume in the proof that each $U_i$ is disjoint from the singular set of $\Lambda$.

Let $D \in C_0$ be a cell with loose chart $U_D$, which is disjoint from $\Lambda \sm D$.  Each loose chart $U_i$ contains two disjointly embedded loose charts \cite{loose}, thus we can assume that there is a Darboux ball $\wt U_D$ which is disjoint from $\Lambda \sm D$, $U_D$ and all $U_i$, and intersects $D$ in the standard plane. We define $\wt D$ to be the Legendrian which is equal to $D$ outside $\wt U_D$, and inside $\wt U_D$ $\wt D$ is equal to a loose Legendrian plane which is formally isotopic to the standard plane and equal to it outside a compact set. Then $\wt D$ is formally isotopic to $D$ via a formal isotopy supported in $\wt U_D$, and both $\wt D$ and $D$ have the same loose chart, $U_D$. Thus $D$ is isotopic to $\wt D$, via an isotopy which fixed on $\Lambda \sm D$. Furthermore, $\wt D$ admits a loose chart which is disjoint from $\Lambda \sm D$ and also disjoint from all $U_i$, namely the loose chart contained in $\wt U_D$. 

We now discard the notation involving the tildes and assume that $D$ admits a loose chart $U_D$ which is disjoint from $\Lambda \sm D$ and also from each $U_i$, the argument above shows that this is no loss of generality.

Because $U_i \sse Y$ is a ball with smooth boundary, there is no obstruction to finding a smooth isotopy, compactly supported on the interior of $D \sm U_D$, which disjoins $D$ from all $U_i$. Using the Serre fibration property for formal Legendrian embeddings over smooth embeddings we see that there is a formal Legendrian isotopy of $D$ with the same property. Since $D$ is loose, Theorem \ref{thm:close loose} implies that there is a Legendrian isotopy $\p_t: D \to Y$ which is compactly supported on the interior of $D$, and which is $C^0$ close to the formal Legendrian isotopy outside of $U_D$. In particular, $\p_1(D)$ is disjoint from all $U_i$. Then $U_i' = \p_1^{-1}(U_i)$ is a loose chart for $D_i$ which is disjoint from $D$, and also disjoint from any cells of $\Lambda$ which $U_i$ was already disjoint from.

We then apply this argument iteratively to each cell in $C_0$, resulting in loose charts for each $D_i$ which are disjoint from all cells in $C_0$.
\end{proof}

\begin{proof}[Proof of Proposition \ref{prop:loose}:]
Let $W_0 \sse \Mor(Q)$ be the union of $W$ and the identities. We define $W_{i+1} \sse \Mor(Q)$ as follows: for any composition $a \overset{f}{\to} b \overset{g}{\to} c \overset{h}{\to} d$ in $Q$ so that $gf, hg \in W_i$, then $f, g, h, gf, hg, hgf \in W_{i+1}$. Then $W_{i+1} \supseteq W_i$, and whenever $W_{i+1} = W_i$ we see that $W_i = \ol W$. Since $\Mor(Q)$ is finite this happens for some finite $i$. Let $C_i \sse C(\Lambda_{A_{n+1}})$ consist of all top-dimensional cells so that $D \in C_i$ exactly when the corresponding morphism $f_D \in \Mor(Q)$ is contained in $W_i$. We prove by induction that all cells in $C_i$ are loose rel $\Lambda_{Q[W^{-1}]}$.

For $i=0$ this is just Proposition \ref{prop:freebloose}. Suppose $D \sse C_{i+1}$, then Lemma \ref{lem:A4 slice} implies there is a set $U \sse \R^n$ and a contactomorphism $(\pi^{-1}(U), \pi^{-1}(U) \cap \Lambda_{A_{n+1}}) \cong (\R^3_\std \x T^*D^{n-2}, \Lambda_6 \x D^{n-2})$, so that $D$ is the cell corresponding to one of the the morphisms $f, g, h$, or $hgf$ in the composition $a \overset{f}{\to} b \overset{g}{\to} c \overset{h}{\to} d$, and the cells corresponding to $gf$ and $hg$ are contained in $C_i$. Lemma \ref{lem:loose slice} implies that $D$ has a loose chart in $\pi^{-1}(U)$, which possibly intersects the cells corresponding to $gf$ and $hg$ but no others. The induction hypothesis is that these cells are themselves loose, and so Lemma \ref{lem:loose induct} implies that $D$ is loose rel $\Lambda_{Q[W^{-1}]}$.
\end{proof}


\begin{thebibliography}{}












\bibitem{affine} R. Casals and E. Murphy, \emph{Legendrian fronts for affine varieties}, arXiv:1610.06977





\bibitem{CE book} K. Cieliebak and Y. Eliashberg, \emph{From Stein to Weinstein and Back: Symplectic Geometry of Affine Complex Manifolds}, Colloquium Publications, 59. AMS, 2012.








\bibitem{EES} T. Ekholm, J. Etnyre, and M. Sullivan, \emph{Non-isotopic Legendrian submanifolds in $\R^{2n+1}$}, J. Differential Geometry, {\bf{71}} (2005), 85--128.




\bibitem{Stein} Y. Eliashberg, \emph{Topological characterization of Stein manifolds of dimension $> 2$}, Internat. J. Math., {\bf{1}} (1990), 29--46.




\bibitem{hbook} Y. Eliashberg and N. Mishachev, \emph{Introduction to the h-Principle}, Graduate Studies in Mathematcs, 48. AMS, 2002.







\bibitem{GKS} S. Guillermou, M. Kashiwara, P. Schapira, \emph{Sheaf quantization of Hamiltonian isotopies and applications to non displaceability problems},  Duke Math. J. 161 (2012), no. 2, 201–-245. 




\bibitem{KS} M. Kashiwara, P. Schapira, {\bf Categories and sheaves}, Grundlehren der Mathematischen Wissenschaften, 332. Springer-Verlag, Berlin, 2006.





\bibitem{loose} E. Murphy, {\em Loose Legendrian embeddings in high dimensional contact manifolds}, preprint, arXiv:1201.2245


\bibitem{N1} D. Nadler, \emph{Arboreal singularities}, Geom. Topol. 21 (2017), no. 2, 1231–-1274. 

\bibitem{N2} D. Nadler, \emph{Non-characteristic expansions of Legendrian singularities}, preprint, arXiv:1507.01513


\bibitem{STZ} V. Shende, D. Treumann, E. Zaslow, \emph{Legendrian knots and constructible sheaves}, Invent. Math. 207 (2017), no. 3, 1031-–1133.

\bibitem{Sta} L. Starkston, \emph{Arboreal singularities in Weinstein skeleta} Selecta Math. (N.S.) 24 (2018), no. 5, 4105-–4140.



\end{thebibliography}
\end{document}